\documentclass[10pt]{amsart} 
\usepackage{amssymb}

\usepackage{graphicx}
\newtheorem{theorem}{Theorem}[section]
\newtheorem{lemma}[theorem]{Lemma}

\theoremstyle{definition}

\theoremstyle{remark}
\newtheorem{remark}[theorem]{Remark}

\theoremstyle{conjecture}
\newtheorem{conjecture}[theorem]{Conjecture}

\numberwithin{equation}{section}

\graphicspath{%
    {converted_graphics/}
    {./}
}

\begin{document}

\title{Towards a statistical proof of the Riemann Hypothesis}


\author{Jon Breslaw}
\address{Department of Economics, Concordia University,
Montreal, Quebec, Canada H3G 1M8}
\email{breslaw@econotron.com}


\date{March 23, 2009.}


\begin{abstract}
Using the $\zeta$ functional equation and the Hadamard product, an analytical expression for the sum of the reciprocal  of the $\zeta$ zeros is established.  We then demonstrate that
on the critical line,  $\left| \zeta \right|$ is convex, and that in the region $0 < \Re(s) \leq 0.5$, 
$\left| \zeta \right|$  has a negative slope, given the RH.  In each case, analytical formulae are established, and numerical examples are presented to validate these formulae. 
\end{abstract}

\maketitle

\section{Introduction}

The Zeta function was first introduced by Euler\ [4] and is defined by: 
\begin{equation}
\zeta (s)=\sum\limits_{n=1}^{\infty }\frac{1}{n^{s}}
\end{equation}

Riemann\ [6] extended this function to the complex plane, for $\Re (s)>0,$
meromorphic on all of $C$, and analytic except at the point $s=1$ which
corresponds to a simple pole: 
\begin{equation}
\zeta (s)=\frac{1}{(1-2^{(1-s)})}\sum\limits_{n=1}^{\infty }\frac{-1^{(n-1)}%
}{n^{s}}
\end{equation}

The Zeta function satisfies the functional equation: 
\begin{equation}
\zeta (s)=\chi (s)\zeta (1-s)
\label{eq:funeq}
\end{equation}
where: 
\begin{equation}
\chi (s)=2^{s}\pi ^{s-1}\sin (.5\pi s)\Gamma (1-s)
\label{eq:chi}
\end{equation}

Let $s^{\ast }$ be a zero root of $\zeta (s)$ such that $\left| \zeta
(s^{\ast })\right| =0$. There are a number of ``trivial zeros'' at $%
s=-2,-4....$ \ All the other (non-trivial) zeros must lie in the critical
strip $0<\Re (s)<1.$\ Clearly, equation \ref{eq:funeq} is satisfied if $\Re (s^{\ast })=.5$.
\ This is the Riemann Hypothesis - the nontrivial zeros of $\zeta (s)$ have
a real part equal to $0.5$. \ The RH is proved if one can demonstrate that
there are no other zero roots - that is, all the roots fall on the critical
line.


 The approach taken in this paper is to first derive analytical expressions for the derivatives of $\left| \chi \right|$; this occurs in Section 3.  Then in Section 4, we derive analytical expressions for the derivatives of  $\left| \zeta \right|$ by exploiting the functional equation and the logarithmic derivative of the Euler product using the Hadamard product. An analytic expression for the sum of the reciprocal of the $\zeta$ zeros is established, which permits a statistically based conjecture about the truth of the RH. We show that on the critical line,  $\left| \zeta \right|$ is convex and that in the region $0 < \Re(s) \leq 0.5$ and $t>2\pi + \Delta $,  $\left| \zeta \right|$  has a negative slope, given the RH.  
 At each stage we demonstrate the validity of the expressions derived using numerical examples.

\section{Analysis}

We follow Riemann's notation: $s = \sigma + it$.  In this work, we take $\Im(s) = t$ as fixed,
following Saidak [7,8].

\begin{lemma}
\label{lemma:zeta0}
\emph{\ }The non trivial zeros of the $\zeta$ function lie on the critical line, or occur in pairs equidistant from the ciritical line.
\end{lemma}

This follows directly  from the functional equation \ref{eq:funeq} - if $\zeta(s)$  is zero, then so is $\zeta(1-s)$. 
Thus the complex zeros of $\zeta$  either lie on the critical line, or are symmetric about it in the strip $0 < \Re(s) < 1 $
Since such zeros occur in pairs, any violation of the RH can be investigated by considering just one side of the critical line.  In this analysis, we consider the area  $\Re(s) \leq 0.5$.

The following properties of the $\zeta$ function are used throughout. On the critical line $s=.5+it$, $\ $%
\begin{eqnarray*}
\Re \zeta (s) &=&\Re \zeta (1-s) \\
\Im \zeta (s) &=&-\Im \zeta (1-s) \\
\arg \zeta (s) &=&-\arg \zeta (1-s)
\end{eqnarray*}

In addition, since we are dealing with  absolute values $\left|  f(s) \right| $, for $s=\sigma+it$, we note:
\begin{equation*}
\ln f(s) = \ln \left |f(s) \right| + i \arg f(s)
\end{equation*}
and hence:
\begin{equation}
\label{eq:abs1}
\frac{ \left| f(s) \right|'} {\left| f(s) \right|} = \Re \left( \frac{f'(s)}{f(s)} \right)
\end{equation}
where $f'(s)$ implies $\frac {\partial f(s)} {\partial \sigma}$ evaluated at $s$.

We will be dealing with a particular region of the critical strip, which we define as $s \in S: 0 < \Re(s) < 1$ and  $\Im(s)>2\pi +\Delta $.


\section{Analysis - $\left| \chi \right| $ function}

The Zeta function satisfies the functional equation: 
\begin{equation}
\zeta (s)=\chi (s)\zeta (1-s)
\label{eq:chi1}
\end{equation}
and hence:
\begin{equation}
\left| \zeta (s)\right|=\left|\chi (s)\right| \left|\zeta (1-s)\right|
\label{eq:chi2}
\end{equation}
where: 
\begin{equation}
\chi (s)=2^{s}\pi ^{s-1}\sin (.5\pi s)\Gamma (1-s)
\label{eq:chi3}
\end{equation}

 In this section, we derive the properties of $\left| \chi(s)\right| $  which are required in ascertaining the properties of $\left| \zeta(s)\right| $.  In particular, we consider the 
analysis of $\left| \chi(s)\right| _{t=t_{0}}$ - the slice of $%
\left| \chi(s)\right| $ when $t$ is fixed at $t_{0}$. 

\ A\ graph of $\left| \chi(s)\right| ,\ s\in S$ is shown in Figure 1 -
the second contour depicts the critical line $\Re (s)=.5$.

\begin{figure}[h] 
  \centering
  \includegraphics[bb=0 0 482 482,width=4.1in,height=4.1in,keepaspectratio]{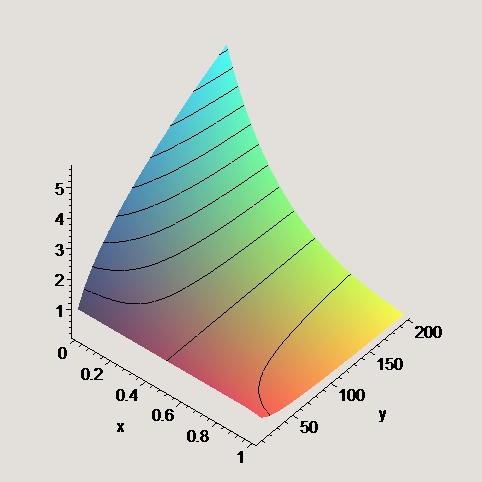}
  \caption{3D\ plot of $\left| \chi(s)\right|$, $\Re(s) = 0...1, \Im(s) = 2\pi..200.$}
  \label{fig:abs2f}
\end{figure}


The following properties of $\left| \chi(s)\right| $ are established:


\begin{lemma}\label{lemma:chi1}
For $s\in S$, $\left| \chi(s)\right| $ is continuous. 
\end{lemma}

From equation \ref{eq:chi3}, for $s\in S,$ $\chi(s)$ is a holomorphic function, and thus $%
\left| \chi(s)\right| $ is continuous.


\begin{lemma}\label{lemma:chi2}
 For $s\in S$, $\frac{\partial }
{\partial \sigma }\left| \chi(s)\right| _{t=t_{0}}=\left| \chi(s)\right| \ \left(
\ln (2\pi )-\Re (\Psi (1-s))\right).$
\end{lemma}

For any $t=t_{0}$: 
\begin{eqnarray*}
\frac{\partial }{\partial \sigma }\left| 2^{s}\right| &=&2^{\sigma }\ln (2)
\\
\frac{\partial }{\partial \sigma }\left| \pi ^{s-1}\right| &=&\pi ^{\sigma
-1}\ln (\pi ) \\
\frac{\partial }{\partial \sigma }\left| \sin (.5\pi s)\right| &=&\frac{%
.5\pi \sin (.5\pi \sigma )\cos (.5\pi \sigma )\ (\cosh ^{2}(.5\pi
t_{0})-\sinh ^{2}(.5\pi t_{0}))}{~\left| \sin (.5\pi s)\right| } \\
&=&\frac{.25\pi \sin (\pi \sigma )}{\left| \sin (.5\pi s)\right| } \\
\frac{\partial }{\partial \sigma }\left| \Gamma (1-s)\right| &=&-\left|
\Gamma (1-s)\right| \ \Re (\Psi (1-s))
\end{eqnarray*}

where $\Psi$ is the digamma function, and the last equation follows from \ref{eq:abs1}.

Taking the modulus of equation \ref{eq:chi3}, we have:
\begin{equation}
\label{eq:chi4}
\left| \chi(s)\right| =\left| \chi(\sigma +it_{0})\right| =\left| 2^{s}\right|
\left| \pi ^{s-1}\right| \left| \sin (.5\pi s)\right| \left| \Gamma
(1-s)\right|
\end{equation}

Hence the derivative of $\left| \chi(s)\right| $ with respect to $%
\sigma $, while keeping $t=t_{0}$, is:
\begin{eqnarray}
\frac{\partial }{\partial \sigma }\left| \chi(s)\right| _{t=t_{0}} &=&\left|
\chi(s)\right| \ (\ \ln (2\pi )+\frac{.25\pi \sin (\pi \sigma )}{\left| \sin
(.5\pi s)\right| ^{2}}-\Re (\Psi (1-s))\ )  \notag \\
&=&\left| \chi(s)\right| \phi (s)
\label{eq:chi5}
\end{eqnarray}

\ The first term in $\phi (s)$, $\ln (2\pi ),$ equals $1.837877067.$ \ The
second term is positive in the critical strip, but quickly becomes
insignificant - \ at $t=2\pi $ this term is less than $0.0001$, and by the
first root $(s=.5+i14.134725)$ it is less than $4\ 10^{-10}$. The last term $%
\Re (\Psi (1-s))$ is the real part of the digamma function (which is
asymptotic with the natural logarithm function) and is positive in the
critical strip for $t>1.1$. For $\Delta =0.02,\ \phi (s)<0$ for $t>2\pi
+\Delta .$ Thus, for $s\in S,$ and  $t>2\pi+\Delta$:
\begin{equation}
\label{eq:chi6}
\frac{\partial }{\partial \sigma }\left| \chi(s)\right| _{t=t_{0}}=\left|
\chi(s)\right| \ \left( \ln (2\pi )-\Re (\Psi (1-s))\right)
\end{equation}
\qquad\ Equation \ref{eq:chi6} works extremely well. \ Using Maple, a numerical
comparison between equation \ref{eq:chi6} and the numerical derivative of $\left| \chi(s)\right| $
evaluated over the critical strip at the first root 
showed a percentage deviation of less than $1.4\ 10^{-7},$ with even lower
deviations for $\sigma \neq .5$. \


\begin{lemma}
\label{lemma:chi3}
 For $s\in S$, the minimum and
maximum of \ $\left| \chi(s)\right| _{t=t_{0}}$\ occurs at $\sigma=1$\ and $\sigma=0$\
respectively.
\end{lemma}

For $s\in S$, $\chi (s)$ in equation \ref{eq:chi5} is a holomorphic function with no zeros.
Hence, the extreme points of $\left| \chi(s)\right| _{t=t_{0}}$will lie on the
boundary.


\begin{lemma}
\label{lemma:chi4}
 For $s\in S$, $\left| \chi(s)\right| > 0.$
\end{lemma}

For $s=\sigma +it_{0},$ let: 
\begin{equation*}
\Phi (s)=\sqrt{(\sin ^{2}(.5\pi \sigma )\cosh ^{2}(.5\pi t_{0})+\cos
^{2}(.5\pi \sigma )\sinh ^{2}(.5\pi t_{0}))}
\end{equation*}
and thus: 
\begin{equation*}
\left| \chi(\sigma +it_{0})\right| =2^{\sigma }\pi ^{\sigma -1}\left| \Gamma
(1-s)\right| \Phi (s)
\end{equation*}
and using standard results\thinspace \lbrack 1]: 
\begin{eqnarray*}
\left| \Gamma (-it_{0}))\right| &=&\sqrt{\frac{\pi t_{0}}{\sinh (\pi t_{0})}}
\\
\left| \Gamma (1-it_{0}))\right| &=&\sqrt{\frac{\pi }{t_{0}\sinh (\pi t_{0})}%
}
\end{eqnarray*}
Hence for $s\in S,$ and from Lemma \ref{lemma:chi3}, 
\begin{eqnarray}
\lim\limits_{t_{0}\rightarrow \infty }\left| \chi(1+it_{0})\right|
&<&\lim\limits_{t_{0}\rightarrow \infty }\left| \chi(\sigma +it_{0})\right|
<\lim\limits_{t_{0}\rightarrow \infty }\left| \chi(0+it_{0})\right|  \notag \\
\sqrt{\frac{2\pi }{t_{0}}} &<&\lim\limits_{t_{0}\rightarrow \infty }\left|
\chi(\sigma +it_{0})\right| <\sqrt{\frac{t_{0}}{2\pi }}
\end{eqnarray}

Thus, $\left| \chi(s)\right| >0.$\ \ These limits  apply even for
relatively low values of $t_{0}$. \ Using Maple, a numerical comparison
between\ $\sqrt{t_{0}/2\pi }$ and $\left| \chi(0+it_{0})\right| $ evaluated at
the first root  showed a percentage deviation of less
than $0.7\ 10^{-7}$.


\begin{lemma}
\label{lemma:chi5}
 For $s\in S$, $t>2\pi +\Delta $, $\frac{\partial }{\partial \sigma }\left| \chi(s)\right| _{t=t_{0}}$ is negative.
\end{lemma}

From Lemma \ref{lemma:chi2}, 
\begin{equation*}
\frac{\partial }{\partial \sigma }\left| \chi(s)\right| _{t=t_{0}}=\left|
\chi(s)\right| \ \left( \ln (2\pi )-\Re (\Psi (1-s))\right)
\end{equation*}
From Lemma \ref{lemma:chi4}, $\left| \chi(s)\right| \ >0$. \ Thus, for small values of $t$ ($%
t<2\pi $), $\frac{\partial }{\partial \sigma }\left| \chi(s)\right| $ is
positive, dominated by the first term, $\ln (2\pi ).$ \ However, for $%
0<\sigma <1,\ t>2\pi +\Delta ,$ the last term $\Re (\Psi (1-s))$ dominates,
and hence $\frac{\partial }{\partial \sigma }\left| \chi(s)\right| _{t=t_{0}}<0$%
. \ This also establishes monotonicity.


\begin{lemma}
\label{lemma:chi6}
If $\Re (s)=.5$, then $\left|\chi(s)\right| =1$.
\end{lemma}

On the critical line, $\left| \zeta (s)\right| =\left| \zeta (1-s)\right| $.
\ From the functional equation \ \ref{eq:funeq}: 
\begin{eqnarray}
\left| \zeta (s)\right| &=&\left| \chi (s)\right| \;\left| \zeta (1-s)\right|
\notag \\
&=&\left| \chi (s)\right| \;\left| \zeta (s)\right|
\end{eqnarray}
 and hence $\left| \chi (s)\right| =1$.

This result also follows from the definition of $\chi (s)$. 
For any $t=t_{0\text{ }}$ and for $\sigma =0.5:$%
\begin{eqnarray*}
\left| 2^{s}\right| &=&2^{\sigma }=\sqrt{2} \\
\left| \pi ^{s-1}\right| &=&\pi ^{\sigma -1}=1/\sqrt{\pi } \\
\left| \sin (.5\pi s)\right| &=&\sqrt{\sin ^{2}(.5\pi \sigma )\cosh
^{2}(.5\pi t_{0})+\cos ^{2}(.5\pi \sigma )\sinh ^{2}(.5\pi t_{0})} \\
&=&\sqrt{\sin ^{2}(\pi /4)\cosh ^{2}(.5\pi t_{0})+\cos ^{2}(\pi /4)\sinh
^{2}(.5\pi t_{0})} \\
\left| \Gamma (1-s)\right| &=&\left| \Gamma (.5-it_{0})\right| =\sqrt{\pi
/\cosh (\pi t_{0})}
\end{eqnarray*}

Thus from equation \ref{eq:chi3}: 
\begin{eqnarray*}
\left| \chi(.5+it_{0})\right| &=&\sqrt{\frac{(2/\pi )(.5\cosh ^{2}(.5\pi
t_{0})+.5\sinh ^{2}(.5\pi t_{0}))\pi }{\cosh (\pi t_{0})}} \\
&=&1 \ 
\end{eqnarray*}

Hence, if $\Re (s)=.5$, then $\left| \chi(s)\right| =1$.


\begin{lemma}
\label{lemma:chi7}
 For $s\in S$, $t>2\pi +\Delta $, $\left| \chi(s) \right|'' = \frac{(\left| \chi(s) \right|')^2}{\left| \chi(s) \right|}$
\end{lemma}

Differentiating  equation \ref{eq:chi5} provides the second derivative of $\left| \chi \right|$:

\begin{equation*}
\left| \chi(s) \right|'' = \left| \chi(s) \right| \left(\frac{-\partial \Re \Psi(1-s)}{\partial \sigma} + \left[\ln(2 \pi) - \Re \Psi(1-s)\right]^2 \right)
\end{equation*}
For even moderate levels of $t$, $\frac {\partial \Re \Psi(1-s)}{\partial \sigma} = 0 $. Hence:
\begin{equation}
\left| \chi(s) \right|'' = \frac{(\left| \chi(s) \right|')^2}{\left| \chi(s) \right|}
\label{eq:chi7}
\end{equation}

Using Maple, a numerical comparison between equation \ref{eq:chi7} and the numerical derivative of equation \ref{eq:chi5}, evaluated at
$\sigma = .5, t=200$,  showed a percentage deviation of less than $1.7\ 10^{-7}$, and even lower off the critical line.

At $\Re (s) = .5$, $\left| \chi(s) \right| = 1$, and thus 

\begin{equation*}
\left| \chi(s) \right|'' = (\left| \chi(s) \right|')^2
\label{eq:chi7b}
\end{equation*}


\begin{lemma}
\label{lemma:chi8}
 For $s\in S$, $t>2\pi +\Delta $, $\left| \chi(s) \right|^{(k)} = \frac{(\left| \chi(s) \right|')^k}{\left| \chi(s) \right|^{k-1}}$
\end{lemma}

The $k$th derivative of $\left| \chi \right|$ can be derived by repeated differentiation of equation \ref{eq:chi7}.  We note that the sign of the $k$th derivative is $(-1)^k$.  Also, for $\Re(s)=.5$, 
$\left| \chi(s) \right|^{(k)} = (\left| \chi(s) \right|')^k$


\section{Analysis - $\left| \zeta \right| $ function}

In this section, we derive the properties of $\left| \zeta(s)\right| $


\begin{lemma}
\label{lemma:zeta1}
 If $s^{\ast }$ is a zero root of $\zeta (s)$, then $\left| \zeta (s^{\ast })\right| =0.$
\end{lemma}

By construction.


\begin{lemma}
\label{lemma:zeta2}
  For $s\in S$, $\Re(s)<0.5,\;\left| \zeta (s)\right| \geq\left| \zeta (\overline{1-s})\right| $
\end{lemma}

 From Lemma \ref{lemma:chi5}, for \emph{\ $s\in S$, }$\left| \chi(s)\right|
_{t=t_{0}}$ is monotonic with $\frac{\partial }{\partial \sigma }\left|
\chi(s)\right| _{t=t_{0}}<0$. \ From Lemma \ref{lemma:chi6} $\left| \chi(0.5+it_{0})\right| =1.
$ Thus for  $\Re (s)<0.5$, $\ \left| \chi(s)\right| >1$.\ Using equation \ref{eq:chi3}, for $\Re
(s)<0.5,$ and noting that $\left| \zeta (s)\right| =\left| \zeta (\overline{s%
})\right| ,$

\begin{eqnarray}
\left| \zeta (s)\right|  &=&\left| \chi(s)\right| \,\left| \zeta (1-s)\right|
=\left| \chi(s)\right| \,\left| \zeta (\overline{1-s})\right|  \\
& \geq &\left| \zeta (\overline{1-s})\right|   \notag
\end{eqnarray}

\begin{remark}
If the strict equality held, we would have RH; however the equality holds at the zeros.
\end{remark}


\begin{lemma}
\label{lemma:zeta3}
For   $\Re (s) < 0.5$, $ \left|\zeta (s)\right|' < 0 $ if  
$ \left|\zeta (1-s)\right|' > 0$.
\end{lemma}

At any point $s$ to the left of the critical line,   $\frac{\partial }{\partial \sigma } \left|\zeta (s)\right|$ can be analytically evaluated as a function of the slope at the corresponding point $1-s$ to the right of the critical line:
\begin{eqnarray}
\zeta (s) & = &\chi (s)\zeta (1-s) \notag \\
\left|\zeta (s)\right| & = &\left|\chi (s)\right| \left|\zeta (1-s)\right| \notag \\
\left|\zeta (s)\right|' & = & -\left|\chi (s)\right| \left|\zeta (1-s)\right|' +
\left|\zeta (1-s)\right|  \left|\chi (s)\right|' 
\label{eq:zeta3}
\end{eqnarray}

$\left|\chi (s)\right|$ and $\left|\zeta (1-s)\right| $ are positive by construction.
By Lemma \ref{lemma:chi5} $\left|\chi (s)\right|' $ is negative always in the critical strip.
Thus if $\left|\zeta (1-s)\right|'$ is positive, then 
$ \left|\zeta (s)\right|'$  has to be  negative. 

\begin{remark}
Unfortunately, this does not help for the situation of interest - that is where 
$\left|\zeta (1-s)\right|' < 0.$
\end{remark}


\begin{lemma}
\label{lemma:zeta4}
  For $\Re(s) = 0.5$,  $\left|\zeta (s)\right| '< 0 $
\end{lemma}

From equation \ref{eq:zeta3},
\begin{equation}
 \left|\zeta (s)\right|  =  
-\left|\chi (s)\right| \left|\zeta (1-s)\right|' +
\left|\zeta (1-s)\right|  \left|\chi (s)\right|' 
\end{equation}

At $s = .5 + it$, $\left| \zeta(s)\right| = \left| \zeta(1-s)\right|$, 
$ \left|\zeta (1-s)\right|' =
  \left|\zeta (s)\right|' $, and by Lemma \ref{lemma:chi6}
$\left|\chi(s)\right| =1$.
Hence, 
\begin{eqnarray*}
\left|\zeta (s)\right|'& = &.5 \left|\zeta (s)\right| \left|\chi (s)\right|' \\
& < &  0
\end{eqnarray*}

since by Lemma \ref{lemma:chi5}
$ \left| \chi(s)\right|' < 0$ 
Thus at $\Re(s) = 0.5$,
\begin{eqnarray}
\label{eq:zeta4}
\frac {\left|\zeta (s)\right|'} {\left|\zeta (s)\right| } = .5 \left|\chi (s)\right|' 
\end{eqnarray}

Using Maple, $\Re \left( \zeta'(s)/\zeta(s)\right) $, evaluated at $s = .5+200i $, gave  -1.7302196, while $.5 \left|\chi (s)\right|' $
gave  -1.7302197.


\begin{lemma}
\label{lemma:zeta5}
$  \sum_{\rho} \Re\left(\frac {1}{\rho} \right)  = -\log(2)-.5\log(\pi)+1+.5\gamma $
where $\rho$ are the non trivial zeros of the $\zeta$ function, and $\gamma$ is Euler's constant. 
\end{lemma}

For any $s$ in ${\mathcal C}$, the $\xi $ function is defined as:
\begin{equation*}
\xi(s) = .5 s (s-1) \pi^{.5 s} \Gamma(.5 s) \zeta(s)
\end{equation*}
The Hadamard product representation of $\xi$ is:
\begin{equation*}
\xi(s) = e^{A+Bs} \prod_{\rho} \left(1-s/\rho\right) e^{s/\rho}
\end{equation*}
where $\rho$  runs over all of the non-trivial zeros of $\zeta(s)$, and A and B are absolute
constants. Taking the logarithmic derivative of $\xi(s)$ we obtain:
\begin{equation*}
\frac{  \xi'(s)} {\xi(s)} = B + \sum_{\rho} \left( \frac{1}{s-\rho}+\frac {1}{\rho} \right)
\end{equation*}
Rewriting in terms of $\zeta$, we obtain:
\begin{equation}
\label{eq:zeta5}
\frac { \zeta'(s)} {\zeta(s)} = B - \frac{1}{s-1} + .5 \log(\pi) -.5\frac{ \Gamma'(.5s + 1)} {\Gamma(.5s+1)} + \sum_{\rho} \left( \frac{1}{s-\rho}+\frac {1}{\rho} \right) 
\end{equation}

Evaluating equation \ref{eq:zeta5} at $s=0$ gives:
\begin{equation*}
B = \log(2)+.5\log(\pi)-1-.5\gamma
\end{equation*}
where $\gamma = -\Gamma'(1)$ is Euler's constant.  Thus $B = -.0230957$.

Since $\left( \frac {\left| \zeta(s) \right|'} {\left| \zeta(s) \right|} \right) = \Re \left(\frac { \zeta'(s)} {\zeta(s)}\right)$,
we can compare equations \ref{eq:zeta4} and \ref{eq:zeta5} evaluated at  $\Re(s) = 0.5$.  The term $\frac{1}{s-1}$ tends to zero as $t$ gets large, and  using Lemma \ref{lemma:chi2}, we note that:
\begin{equation} 
\label{eq:zeta5a}
\ln(\pi)-\Re (\psi(.5s+1)) =\ln(2 \pi)-\Re (\Psi(1-s)) = \frac {\left| \chi \right|'}  {\left| \chi \right|};
\end{equation} 
Thus, if these two equations are to have the same value at $\Re(s) = 0.5$, then:  
\begin{equation} 
\label{eq:zeta5b}
\sum_{\rho} \Re \left( \frac{1}{s-\rho}+\frac {1}{\rho} \right)  = \Re \left(\frac {1}{\rho} \right) = -B
\end{equation} 
Evaluating the sum over the first 100000 non trivial $\zeta$ zeros for $s = .5 +it$ for a range of $t$, $t> 200$, generated a constant value of 0.023073645, which is approximately -B.\footnote{The difference comes about because of the summation over only 100000 terms, instead of infinity. We noted also that the sum was over both positive and negative values of the zeros - and indeed in Hadamard's work the three regions were defined in terms of the absolute value of $\rho$.  In addition, the analysis is not valid for a situation where $s$ is equal to $\rho_k$, a $\zeta$ zero, since $\left| \zeta(s) \right|$ is not analytic at that point.  However it is valid in the limit.}


\begin{lemma}
\label{lemma:zeta6}
  For $\Re(s) = 0.5$,  $\left|\zeta (s)\right|'' > 0 $.
\end{lemma}

Differentiating equation \ref{eq:zeta5}:
\begin{equation}
\label{eq:zeta6}
\frac { \zeta''(s)} {\zeta(s)} - \left(\frac { \zeta'(s)} {\zeta(s)}\right)^2 = \frac{1}{(s-1)^2}  -\frac{1}{8} \Psi'(.5s + 1) - \sum_{\rho} \left( \frac{1}{(s-\rho)^2} \right) 
\end{equation}

Note that $\lim_{t\rightarrow \infty }\frac{1}{(s-1)^2} = 0$ and $\lim_{t\rightarrow \infty }\Psi'(s) = 0$. In addition, evaluated at  $\Re(s) = 0.5$, $ \Re \left( \sum_{\rho} \left[ \frac{1}{(s-\rho)^2} \right] \right) = \Re (K(2,s)) < 0 $.  Thus we have:

\begin{equation}
\frac { \zeta''(s)} {\zeta(s)} = \left(\frac { \zeta'(s)} {\zeta(s)}\right)^2 - K(2,s) > 0
\end{equation}
and equivalently:
\begin{equation}
\label{eq:zeta6b}
\frac { \left|\zeta(s) \right|''} {\left| \zeta(s) \right|} = \left(\frac {\left| \zeta(s)\right|'} {\left|\zeta(s)\right|}\right)^2 - \Re(K(2,s))
\end{equation}

Again, we check to see that this is correct numerically. Using $s = .5 + 200i$, $\Re \left( \frac { \zeta''(s)} {\zeta(s)} - \left(\frac { \zeta'(s)} {\zeta(s)}\right)^2 \right) = 1.4511544$.  Evaluating $ \sum_{\rho} \left( \frac{1}{(s-\rho)^2} \right)$ over 100000 zeros generated  $K(2,s) = -1.4511540$.

\begin{remark}
 The evaluation of $K(2,s) = \sum_{\rho} \Re \left( \frac{1}{(s-\rho)^2} \right) $, is only a function of $t$, since $\Re(s)= \Re(\rho) $. Defining $K(k,s) = \Re \left(\frac {\partial^k}{\partial \sigma^k} \frac {1}{(s-\rho)} \right)$, it is clear that  for $k$ odd and $\Re(s) = .5$, $\Re(K(k,s)) = 0$.  
\end{remark}


\begin{lemma}
\label{lemma:zeta8}
  For $0 < \Re(s) < 0.5$ and  $t>2\pi +\Delta $, and given the RH,  $\left|\zeta (s)\right|' < 0 $.
\end{lemma}

Applying Lemma \ref{lemma:zeta5} and equation \ref{eq:zeta5a} to equation \ref{eq:zeta5}  gives :
\begin{equation*}
 \frac {\left| \zeta(s) \right|'} {\left| \zeta(s) \right|} =  .5 \frac {\left| \chi(s) \right|'}  {\left| \chi(s) \right|}
 + \sum_{\rho} \Re \left( \frac{1}{s-\rho} \right) 
\end{equation*}
where as before we ignore the term  $ \Re \frac{1}{s-1} $ which becomes negligible for large $t$. Rewriting as:
\begin{equation}
\label{eq:zeta8}
 \left| \zeta(s) \right|' =  .5  \left| \zeta(s) \right| \left(  \frac {\left| \chi(s) \right|'}  {\left| \chi(s) \right|}
 + \sum_{\rho} \Re \left( \frac{1}{s-\rho} \right)  \right) 
\end{equation}
demonstrates how equation \ref{eq:zeta4} has been extended to cover the range $0 < \Re(s) < 1$. 

It is easy to check the numerical validity of equation \ref{eq:zeta8} by using the relationship shown in equation \ref{eq:abs1}:
\begin{equation}
\label{eq:zeta8b}
 \left| \zeta(s) \right|' = \left| \zeta(s) \right|  \Re \left( \frac{\zeta'(s)}{\zeta(s)} \right)
\end{equation}
Evaluating  $\left| \zeta(s) \right|'$ at $s = .1 + 200i$ using equation \ref{eq:zeta8b} gives -28.52551836, while evaluating  $\left| \zeta(s) \right|'$  using equation \ref{eq:zeta8} gives  -28.52551645, where the sum over zeros (-0.55098684) is evaluated over the first 100000 zeros. A second example uses $t =9291.071149949$, such that $\zeta(.5+it)$ is close to a zero. 
For $s = .1 + it$, the sum over zeros is -3.8557529, with the two estimates being -199.4723 and -199.4719 respectively.  Flipping to the other side, with $s = .9 + it$, the sum over zeros is reversed (3.8557529), and the two estimates are 0.2958258 and 0.2958254 respectively.

 Consider equation \ref{eq:zeta8}, and recall that $\Re \left(\frac{1}{a+ib}\right) =  \frac{a}{a^2+b^2}$. 
\begin{enumerate}
  \item  On the critical line, equations  \ref{eq:zeta4} and \ref{eq:zeta8} are identical, since  at $\Re(s) = .5$, $\left| \chi(s) \right| = 1$ by Lemma \ref{lemma:chi6}, and $ \Re \left( \frac{1}{s-\rho} \right) = 0 $ since $\Re(s) =\Re(\rho) = .5. $
  \item  The sign of   $\left| \zeta(s) \right|'$  in  equation \ref{eq:zeta8} in the range $0 < \Re(s) < .5$  depends on the sign of the terms in the brackets.  The first term, $\frac {\left| \chi(s) \right|'}  {\left| \chi(s) \right|}$ is always negative, since  $\left| \chi(s) \right|'$ is negative and monotonic by Lemma \ref{lemma:chi5}. The sign of the second term, involving the sum over the $\zeta$ zeros, depends on the validity of the RH. 
  For a single zero, with $\Re(\rho) = .5$, the contribution is negative since in this range, $\Re(s - \rho) $ is negative.
  For a symmetric pair of zeros (supposing  RH were false), say $\rho_L$ and $\rho_R$, the contribution will be  negative only if $s$ is to the left of $\rho_L$.  If $s$ is to the right of $\rho_L$ and to the left of the critical line, that is $\Re(\rho_L) < \Re(s) < .5$,  then the total contribution of the pair is positive.

  Thus the RH cannot be established from a consideration of just the Hadamard product.
\begin{qed}
\hspace{1.5in}
\end{qed}

\end{enumerate}

\begin{conjecture}
The zeros of the $\zeta$ function are derived from a single stochastic process.
\end{conjecture}
  Since $\lim_{\Re(s) \rightarrow \infty}  \zeta(s)  = 1$, it follows directly from the functional equation \ref{eq:funeq} that $\lim_{\Re(s) \rightarrow -\infty}   \zeta(s) = \chi(s)$. Equation \ref{eq:zeta8} demonstrates  that, with the exception of the region around a $\zeta$ zero, there is a high correlation between the behaviour of the $\zeta$ function and the $\chi$ function.  The critical difference between the two is that by Lemma \ref{lemma:chi4},  $ \chi(s) $ has no zeros in $S$. We know that there are no $\zeta$ zeros to the left of $\Re(s) = 0$, so the question is whether the $\zeta$ function has become sufficiently close to the $\chi$ function to the left of  $\Re(s) = .5$ such that there are no zeros to the left of the critical line.

  We can get an indication as to whether the RH is true from a consideration of Lemma \ref{lemma:zeta5}.  We know from the work on random matrix theory that the distribution of $\zeta$ zeros on the imaginary axis is a statistical phenomena.  Thus, should these zeros also occur as  symmetric pairs, then this too will also be a statistical phenomena.  Let the $n$th $\zeta$ zero, $\rho_n = \sigma_n + it_n$, and consider the statistic:
\begin{eqnarray}
 \sum_{\rho} \Re\left(\frac {1}{\rho} \right)& = &\sum_n \Re\left(\frac {\sigma_n}{\sigma_n^2+t_n^2} \right)\\
                                               & \approx &\sum_n \Re\left(\frac {\sigma_n}{t_n^2} \right)
\end{eqnarray}
Hence, if the RH is false, then  $ \sum_{\rho} \Re\left(\frac {1}{\rho} \right)$ consists of terms that are the ratio of two statistically derived variates - $\sigma_n$ and $t_n$. If we consider some standard statistical distributions that are similarly functions of the ratio of stochastic processes - beta (ratio of two gamma processes), F (ratio of two $\chi^2$ processes), Student's t (ratio of normal and $\sqrt(\chi^2)$ processes) - we observe that the distribution function cannot be expressed in an analytic form, although the converse is not necessarily true.

By Lemma \ref{lemma:zeta5} we note that $  \sum_{\rho} \Re\left(\frac {1}{\rho} \right)  = -\log(2)-.5\log(\pi)+1+.5\gamma $
where $\rho$ are the non trivial zeros of the $\zeta$ function, and $\gamma$ is Euler's constant. The fact that this sum can be expressed as an analytic form solely in terms of mathematical constants - and this is of course the crux of the conjecture - would imply that  $\sum_n \Re\left(\frac {\sigma_n}{t_n^2} \right)$ does not consist of two stochastic processes.  That would of course be true only if $\sigma_n = 0.5 \ \forall \ n$ - which would give us the RH.


\section*{Acknowledgment}
I would like to thank Marco Bertola, Ken Davidson, Lubo\u{s} Motl and Terence Tao for their
valuable comments and feedback during the preparation of this article.

\bibliographystyle{amsplain}

\begin{thebibliography}{10}


\bibitem {A} T. Aoki, \textit{Calcul exponentiel des op\'erateurs
microdifferentiels d'ordre infini.} I, Ann. Inst. Fourier (Grenoble)
\textbf{33} (1983), 227--250.


\bibitem {1} M. Abramowitz and I.A. Stegun, (Eds.), (1972). \emph{Handbook of
Mathematical Functions with Formulas, Graphs, and Mathematical Tables}, New
York: Dover, (1972), Eqn 6.1.29-31.

\bibitem {2} J.B. Conrey \emph{The Riemann Hypothesis}, \textbf{50}, \
http://www. ams.org/notices/200303/fea-conrey-web.pdf 

\bibitem {3} L. Euler, \emph{Introductio in Analysin Infinitorum, }Chapter
15, Lausanne. (1748).

\bibitem {4} W. Feller, ``Stirling's Formula.'' \S 2.9 in \emph{An
Introduction to Probability Theory and Its Applications}, \textbf{1}, 3rd ed.
New York: Wiley, (1968), 50-53.

\bibitem {5} G.F.B. Riemann,  ``\"{U}ber die Anzahl der Primzahlen
unter einer gegebenen Gr\"{o}sse.'' \emph{Monatsber}, K\"{o}nigl. Preuss.
Akad. Wiss. Berlin, (1859), 671-680.


\bibitem {6} F. Saidak and P. Zvengbrowski, (2003). ``On the Modulus of the
Riemann Zeta Function in the Critical Strip'',\emph{\ Math Slovaca}, Vol
53(2), 145-172.

\bibitem {7} F. Saidak, (2004). ``On the Logarithmic Derivative of the Euler Product'',
\emph{\ Tatra Mt. Math. Publ.} 113-122.


\end{thebibliography}

\end{document}